\documentclass[12pt]{amsart}

\newcommand{\vv}{v^*_1}
\newcommand{\vvv}{v^*_2}
\newcommand{\vz}{v^*}
\newcommand{\ttt}{\Theta}
\newcommand{\ot}{\otimes}
\newcommand{\vi}{v_i}
\newcommand{\vzi}{v^*_i}
\newcommand{\Io}{I_{\bar 0}}
\newcommand{\IIo}{I_{\bar 1}}
\newcommand{\ei}{e_i}

\newcommand{\eij}{e_{ij}}
\newcommand{\eji}{e_{ji}}
\newcommand{\fij}{f_{ij}}
\newcommand{\fji}{f_{ji}}
\newcommand{\fii}{f_{ii}}

\newcommand{\eii}{e_{ii}}
\newcommand{\ejj}{e_{jj}}
\newcommand{\al}{\alpha}
\newcommand{\ds}{\displaystyle}
\newcommand{\vh}{\varphi}
\newcommand{\vr}{\varepsilon}

\newcommand{\eal}{e^{\frac12\al}}
\newcommand{\ealm}{e^{-\frac12\al}}
\newcommand {\be}{\beta}
\newcommand{\de}{\Delta}
\newcommand{\dekip}{\Delta_i^{(k)+}}

\newcommand{\dekim}{\Delta_i^{(k)-}}

\newcommand{\dekiip}{\Delta_i^{(k-1)+}}
\newcommand{\dekiim}{\Delta_i^{(k-1)-}}
\newcommand{\dekjip}{\Delta_j^{(k-1)+}}
\newcommand{\dekjim}{\Delta_j^{(k-1)-}}

\newcommand {\supplus}{\mathop{{\supset}\llap{\raise 
0.5pt\hbox{\normalfont\small+}\hskip 0.5pt}}} 

\newcommand {\subplus}{\mathop{{\subset}\llap{\raise 
0.5pt\hbox{\normalfont\small+}\hskip 0.5pt}}}  

\newcommand{\pp}{\partial}

\newcommand {\Cee}    {{\mathbb  C}}

\newcommand {\Zee}    {{\mathbb  Z}}

\newcommand {\fA}     {{\mathfrak{A}}}

\newcommand {\fb}     {{\mathfrak{b}}}

\newcommand {\fg}     {{\mathfrak{g}}}    %
\newcommand {\fgl}    {{\mathfrak{gl}}}  %
\newcommand {\fh}     {{\mathfrak{h}}}

\newcommand {\fn}     {{\mathfrak{n}}}

\newcommand {\fpe}    {{\mathfrak{pe}}}   %

\newcommand {\fq}     {{\mathfrak{q}}}

\newcommand {\fS}     {{\mathfrak{S}}}

\newcommand {\cal} {\mathcal}

\def \opname#1#2%
  {\expandafter\newcommand \csname #1\endcsname {{\mathop{#2}\nolimits}}}

\newcommand{\rmname}[1]
  {\expandafter\newcommand \csname #1\endcsname {{\operatorname{#1}}}}

\newcommand{\rmnameii}[2]
  {\expandafter\newcommand \csname #1\endcsname {{\operatorname{#2}}}}

\rmname{act}
\rmname{Ad}
\rmname{Add}
\rmname{ad}
\rmname{Alt}
\rmname{alt}
\rmname{Ann}
\rmname{antidiag}
\rmname{Ber}
\rmname{ber}
\rmname{Br}
\rmname{card}
\rmname{ch}
\rmname{Char}
\rmname{cem}
\rmname{cj}
\rmname{Cliff}
\rmname{cntr}
\rmname{codim}
\rmname{coind}
\rmname{const}
\rmname{col}
\rmname{cork}
\rmname{cpr}
\rmname{diag}
\rmnameii{Div}{div}
\rmname{Def}
\rmname{Der}
\rmname{Dim}
\rmname{End}
\rmname{Even}
\rmname{Ext}
\rmname{gr}
\rmname{Hom}
\rmname{HT}
\rmnameii{Ht}{ht}
\rmname{hwt}
\rmname{Id}
\rmname{id}
\rmname{ind}
\rmname{Ind}
\rmname{Inf}
\rmname{irr}
\rmname{Le}
\rmname{Lie}
\rmname{lwt} 
\rmname{mult}
\rmname{Maps}
\rmname{Mor}
\rmname{nm}
\rmname{Ob}
\rmname{Odd}
\rmname{Osc}
\rmname{per}
\rmname{Pic}
\rmname{pr}
\rmname{pro}
\rmname{Prime}
\rmname{Proj}
\rmname{prt}
\rmname{pt}
\rmname{Q}
\rmname{qet}
\rmname{qtr}
\rmname{rd}
\rmname{rk}
\rmname{row}
\rmname{Res}
\rmname{salt}
\rmname{Sch}
\rmname{SBr}
\rmname{scalar}
\rmname{Ser}
\rmname{sign}
\rmname{Smbl}
\rmname{spin}
\rmname{ssym}
\rmname{str}
\rmname{st}
\rmname{sgn}
\rmname{sq}
\rmname{symm}
\rmname{supp}
\rmname{Supp}
\rmname{St}
\rmname{Spec}
\rmname{Spm}
\rmname{tr}
\rmname{vpt}
\rmname{weyl}
\rmname{Weyl}
\rmname{Witt}

\opname{vvol}  {{v\hspace{-0.1ex}o\hspace{-0.02ex}l\/}}
\opname{pnt}  {\text{\normalfont pt}}
\opname{Span} {{Span}}
\opname{slim} {\overline{\lim}}
\opname{Vol}  {{V\hspace{-0.55ex}o\hspace{-0.02ex}l\/}}
\opname{QVol} {{Q\hspace{-0.3ex}V\hspace{-0.55ex}o\hspace{-0.02ex}l\/}}
\opname{PoVol}{{P\hspace{-0.35ex}o\hspace{-0.25ex}V\hspace{-0.55ex}o\hspace{-0.02ex}l\/}}
\opname{BVol} {{B\hspace{-0.2ex}V\hspace{-0.55ex}o\hspace{-0.02ex}l\/}}
\opname{Par}  {{P\hspace{-0.3ex}a\hspace{-0.05ex}r\/}}

\rmname{Mat}
\rmname{Bil}
\rmname{Diff}
\rmname{Ker}
\rmname{Herm}
\rmname{Coker}
\rmname{Conn}
\rmname{Covect}
\rmname{Vect}
\rmname{Int}

\rmnameii {IM} {Im}
\rmnameii {RE} {Re}

\opname{Aut} {{A\hspace{-0.2ex}u\hspace{-0.1ex}t\/}}
\opname{GL} {{G\hspace{-0.3ex}L}}
\opname{SL} {{S\hspace{-0.3ex}L}}
\opname{Exp} {{E\hspace{-0.2ex}x\hspace{-0.1ex}p\/}}
\opname{GQ} {{G\hspace{-0.2ex}Q}}
\opname{OSp} {{O\hspace{-0.25ex}S\hspace{-0.15ex}p\/}}
\opname{Out} {{O\hspace{-0.25ex}u\hspace{-0.15ex}t\/}}
\opname{Spp} {{S\hspace{-0.2ex}p\/}}
\opname{SpO} {{S\hspace{-0.2ex}p\hspace{-0.02ex}O\/}}
\opname{Pe} {{P\hspace{-0.25ex}e\/}}
\opname{SPe} {{S\hspace{-0.25ex}P\hspace{-0.25ex}e\/}}
\opname{Spin} {{S\hspace{-0.25ex}p\hspace{-0.05ex}i\hspace{-0.1ex}n\/}}
\opname{Iso} {{I\hspace{-0.25ex}s\hspace{-0.1ex}o\/}}
\opname{SSPe} {{S\hspace{-0.25ex}S\hspace{-0.15ex}P\hspace{-0.25ex}e\/}}
\opname{PeU} {{P\hspace{-0.25ex}e\hspace{-0.1ex}U\/}}
\opname{QU} {{Q\hspace{-0.15ex}U\/}}
\opname{U} {{U\/}}

\opname{cGQ} {{\cal G \hspace{-0.2em} Q \/}}
\opname{cSL} {{\cal S \hspace{-0.2em} L \/}}
\opname{cGL} {{\cal G \hspace{-0.2em} L \/}}
\opname{cOSp} {{\cal O \hspace{-0.2em} S \hspace{-0.3em} \it p\/}}
\opname{cPe} {{\cal P \hspace{-1.5pt} \it e\/}}
\opname{cVect} {{\cal V \hspace{-1.5pt} \it e\hspace{-0.1ex}c\hspace{-0.1ex}t\/}}
\opname{cVol} {{\cal V \hspace{-1.5pt} \it o\hspace{-0.1ex}l\/}}
\opname{cAut} {{\cal A \hspace{-0.2em} \it u\hspace{-0.1em}t\/}}
\opname{cCovect} {{\cal C \hspace{-1.5pt}
     \it o\hspace{-0.1ex}v\hspace{-0.1ex}e\hspace{-0.1ex}c\hspace{-0.1ex}t\/}}
\opname{CW} {{C\hspace{-0.15ex}W}}

\opname {Ab}   {{\sf Ab}}
\opname {Alg}   {{\sf Alg}}
\opname {ASch}  {{\sf Aff\;Sch}}
\opname {Funct}   {{\sf Funct}}
\opname {Gr}   {{\sf Gr}}
\opname {Grf}  {{{\sf Gr}_f}}
\opname {Mods}   {{\sf mods}}
\opname {Rings}   {{\sf Rings}}
\opname {Salg}   {{\sf Salg}}
\opname {Sets} {{\sf Sets}}
\opname {SSMan} {{\sf SMan}}
\opname {Top}  {{\sf Top}}
\opname {Vebun}   {{\sf Vebun}}

\newcommand {\ev} {{\bar0}}

\newcommand {\eps} {\varepsilon}

\newcommand {\tto} {\longrightarrow}
\newcommand {\pder}[1] {{\frac{\partial}{\partial {#1}}}}

\newcommand {\bcdot}   {\mathbin{\hbox{\raise.4ex\hbox{\bf.}}}} 

\newcommand {\secno} {}
\newcommand {\ssecfont} {\normalfont\bf}

\newtheorem{Theorem}{\secno Theorem}

\newtheorem{Lemma}[Theorem]{\secno Lemma}

\newenvironment {th*}[1]
    {\gdef\thname{#1} \begin{thn}}%
    {\end{thn}}
\newtheorem{thn}[Theorem] {\thname}

\theoremstyle{definition}

\newtheorem{Convention}[Theorem]{\secno Convention}

\newenvironment {ex*}[1]
    {\gdef\thname{#1} \begin{exn}}%
    {\end{exn}}
\newtheorem{exn}[Theorem]{\thname}

\theoremstyle{remark}

\newenvironment {rem*}[1]
    {\gdef\thname{#1} \begin{remn}}%
    {\end{remn}}
\newtheorem{remn}[Theorem]{\thname}

\newcommand {\ssec}{\subsection*}

\newcommand {\ssbegin}[2]
  {\def \secno {\gdef \secno {}{\ssecfont #1. }}%
   \begin{#2}}
\begin{document}
    
\title[Projective Schur functions and superspaces]
{Projective Schur functions as a bispherical functions on certain 
homogeneous superspaces}

\author{Alexander Sergeev} 

\address{ Balakovo Institute of Technique, Technology and Control,Chapaeva 140
, Balakovo, Saratov Region, Russia}
\email{ sergeev@bittu.org.ru}

\thanks{I am thankful to G.~Olshanskii who raised the problem ,
D.~Leites for help and Isaac Newton Institute for hospitality and support}



\begin{abstract} I show that the projective Schur functions may be
interpreted as bispherical functions of either the triple $(\fq(n), 
\fq(n)\oplus \fq(n), \fq(n))$, where $\fq(n)$ is the ``odd'' (queer)
analog of the general linear Lie algebra, or the triple $(\fpe(n), 
\fgl(n|n), \fpe(n))$, where $\fpe(n)$ is the periplectic Lie
superalgebra which preserves the nondegenerate odd bilinear form
(either symmetric or skew-symmetric). Making use of this
interpretation I characterize projective Schur functions as common
eigenfunctions of an algebra of differential operators.
\end{abstract}

\maketitle

\section*{Introduction}

The ground field is $\Cee$.

\ssec{1.1} In \cite{Sch} I.~Schur introduced projective Schur
functions as characteristics of projective representations of
symmetric groups. In \cite{Se1} I showed that these characteristics
are actually characters of tensor representations of Lie superalgebra
$\fq(n)$.

In \cite{St} Stembridge interpreted projective Schur functions as
characteristics of of spherical functions of a certain twisted Gelfand
pair (for definition see \cite{St}).  Here I will remind definition of
bispherical functions and demonstrate that these Stembridge's
characteristics are precisely bispherical functions of the triple
$(\fpe(n), \fgl(n|n), \fpe(n))$, where $\fpe(n)$ can be embedded into
$\fgl(n|n)$ in two different ways (corresponding to interpretation of
$\fpe(n)$ as the algebra preserving a symmetric odd bilinear form or a
skew-symmetric one).

One obtains one more realization of projective Schur functions if one
considers bispherical functions of the triple $(\fq(n), \fq(n)\oplus
\fq(n), \fq(n))$, where $\fq(n)$ can be embedded into $\fq(n)\oplus
\fq(n)$ by one of the two ways: either as the diagonal or as a
``twisted diagonal''.

Both ways to realize projective Schur functions allow one to construct
an algebra of differential operators for which the projective Schur
functions are eigenfunctions.  This algebra appears as the algebra of
radial parts of Laplace operators for the Lie superalgebra
$\fgl(n|n)$ or $\fq(n)$.  

\ssec{1.2. Differential operators and projective Schur functions} Let
$I=\{1, \dots, n\}$, $V$ an $n$-dimensional vector space, 
$\{e_i\}_{i\in I}$ a basis of $V$, and $\{\eps_i\}$ the dual basis of
$V^*$. If $l\in V^*$, than $e^l$ denotes a homomorphism $S(V)\tto
\Cee$, where $S(V)$ is the symmetric
algebra of $V$. Recall that $S(V)^*$ can be identified with the
algebra of formal power series in $n$ indeterminates.

On $S(V)^*$, define a family of differential operators. Set
$\vr_{ij}=\vr_i-\vr_j$ for $i, j\in I, \ i\ne j$ and set
$$
\pp_i^{(1)}(e^l)=l(e_i)e^l=\pp_i(e^l) \text{ for any } l\in V^*, 
$$
set further
$$
\partial_i^{(k)}=\begin{cases}\partial_i\partial_i^{(k-1)}
+\sum\limits_{j\ne i}\displaystyle\frac{2}
{e^{\varepsilon_{ij}}-e^{\varepsilon_{ji}}}
(\partial_i^{(k-1)}-\partial_j^{(k-1)})&\text{ for $k$ odd}, \cr
(\partial_i-1)\partial_i^{(k-1)}
+\sum\limits_{j\ne i}\left(\displaystyle\frac{2}
{e^{\varepsilon_{ij}}-e^{\varepsilon_{ji}}}
\partial_i^{(k-1)}-\displaystyle
\frac{2e^{\varepsilon_{ij}}}{e^{\varepsilon_{ij}}
-e^{\varepsilon_{ji}}}
\partial_i^{(k-1)}\right)&
 \text{ for $k$ even}.\end{cases}\eqno{(1.2.1)}
$$
Finally, for $k$ odd, set 
$$
\Omega_k=\sum\limits_{i=1}^n\partial_i^{(k)}.
$$
It is not difficult to verify directly that 
$$
\renewcommand{\arraystretch}{1.4}
\begin{array}{l}
\Omega_3=\ds\sum\limits_1^n\partial_i^3+\ds\sum\limits_{i<j}
\ds\frac{6}
{e^{\varepsilon_{ij}}-e^{\varepsilon_{ji}}}(\partial_i^2-\partial_j^2)
-\ds\sum
\ds\frac{6}
{(e^{\frac12\varepsilon_{ij}}+
e^{\frac12\varepsilon_{ji}})^2}(\partial_i+\partial_j)\\
+24\ds\sum\limits_{i\notin \{j, k\}}
\ds\frac{1}{(e^{\varepsilon_{ij}}-e^{\varepsilon_{ji}})
(e^{\varepsilon_{ik}}-e^{\varepsilon_{ki}})}
\partial_i-\left(\sum\limits_1^n\partial_i\right)^2, 
\end{array}
\eqno(1.2.2)
$$
where $\{j, k\}\subset I$ is any two-element subset. 

Introduce new indeterminates: $x_i=e^{\varepsilon_i}$. Then
$\partial_{i}= x_i\pder{x_i}$ and
$$
\partial_i^{(k)}=\begin{cases}\partial_i\partial_i^{(k-1)}
+\sum\limits_{j\ne
i}\displaystyle\frac{2x_ix_j}{x_i^2-x_j^2}\left(\partial_i^{(k-1)}-
\partial_j^{(k-1)}\right)& \text{ for $k$ odd}, \cr
(\pp_i-1)\pp_i^{(k-1)} +\sum\limits_{j\ne
i}\left(\displaystyle\frac{2x_ix_j}{x_i^2-x_j^2}\pp_i^{(k-1)}
-\displaystyle\frac{2x_i^2}{x_i^2-x_j^2}\pp_j^{(k-1)}\right) &\text{ for $k$
even}.\end{cases}\eqno{(1.2.3)}
$$
We have
$$
\renewcommand{\arraystretch}{1.4}
\begin{array}{l}
\Omega_3=\ds\sum\pp_i^3+6\ds\sum\limits_{i<j}\ds
\frac{x_ix_j}{x_i^2-x_j^2}(\pp_i^2-\pp_j^2)-6\ds\sum\limits_{i<j}\ds\frac{x_ix_j}{(x_i+x_j)^2}(\pp_i+\pp_j)
\\
+24\ds\sum\limits_{i<j<k}\ds\frac{x_ix_jx_k}{(x_i^2-x_j^2)
(x_i^2-x_k^2)}\left 
(x_i(x_j^2-x_k^2)\pp_i-x_j(x_i^2-x_k^2)\pp_j+x_k(x_i^2-x_j^2)\pp_k
\right)
\\
-\left(\ds\sum\partial_i\right)^2.
\end{array}
\eqno(1.2.4)
$$

For $i=1, \dots , n$ and any $k$, define differential operators
$\tilde{\pp}_i^{(k)}$ and $\tilde{\pp}_{\bar \imath}^{(k)}$ by setting
$$
\tilde {\pp}_i^{(1)}= \tilde {\pp_i}=\tilde {\pp}_{\bar \imath}^{(1)}
= \tilde {\pp}_{\bar \imath}=\pp_i
$$ 
and
$$
\renewcommand{\arraystretch}{1.4}
\begin{array}{rcl}
\tilde {\pp}_i^{(k)}&=& \tilde {\pp}_i \tilde {\pp}_i^{(k-1)}+
\ds\sum\limits_{j\ne i} \ds\frac{e^{\frac12\varepsilon_{ij}}}{
e^{\frac12\varepsilon_{ij}}- e^{\frac12\varepsilon_{ji}}}
(\tilde{\pp}_i^{(k-1)}-\tilde{\pp}_j^{(k-1)}) \\
&-&\ds\sum\limits_{j\ne i}\ds\frac{e^{\frac12\varepsilon_{ij}}}{
e^{\frac12\varepsilon_{ij}}+ e^{\frac12\varepsilon_{ji}}}
(\tilde{\pp}_i^{(k-1)}+\tilde{\pp}_{\bar \jmath}^{(k-1)}),
\end{array}
\eqno(1.2.5)
$$
$$
\renewcommand{\arraystretch}{1.4}
\begin{array}{rcl}
\tilde {\pp}_{\bar \imath}^{(k)}&=& -\tilde {\pp}_i \tilde {\pp}_{\bar
\imath}^{(k-1)} -\ds\sum\limits_{j\ne
i}\ds\frac{e^{\frac12\varepsilon_{ij}}}{ e^{\frac12\varepsilon_{ij}}-
e^{\frac12\varepsilon_{ji}}} (\tilde{\pp}_{\bar
\imath}^{(k-1)}-\tilde{\pp}_{\bar \jmath}^{(k-1)})\\
&+&\ds\sum\frac{e^{\frac12\varepsilon_{ij}}}{
e^{\frac12\varepsilon_{ij}}+ e^{\frac12\varepsilon_{ji}}}
(\tilde{\pp}_{\bar \imath}^{(k-1)}+\tilde{\pp}_j^{(k-1)}).
\end{array}
\eqno(1.2.6)
$$
Set further
$$
\widetilde\Omega_{k}=
\sum\limits_{i=1}^n(\tilde{\pp}_i^{(k)}+ \tilde{\pp}_{\bar \imath}^{(k)}).
$$

\ssbegin{1.2.1}{Lemma} The algebra $\Omega$ generated by operators $\Omega_k$
for $k=1, 3, 5, \dots$ coincides with the algebra generated by
operators $\widetilde\Omega_{k}$ for $k=1, 2, 3, 4, \dots$. 
\end{Lemma}

\ssec{1.2.2} Let $x_1, \ldots , x_n, t$ be indeterminates. 
Introduce polynomials $q_k(x_1, \ldots , x_n)$ from equation
$$
\sum\limits_{k=0}^{\infty} q_k(x_1, \ldots , x_n)t^k
=\frac{\prod(1+x_it)}{\prod(1-x_it)}.\eqno(1.2.7)
$$
Further, set 
$$
Q_{k, l}=q_kq_l+2\sum\limits_{p=1}^{l} q_{k+p}q_{l-p}. 
$$
We see that $Q_{(k, 0)}=q_k$ and $Q_{k, l}= -Q_{l, k}$ for $k+l>0$. If
$\lambda=(\lambda_1, \ldots , \lambda_m)$ is a strict partition, i.e.,
$\lambda_1>\lambda_2>\ldots>\lambda_m>0$, then ($Pf$ is the Pfaffian 
of the skew-symmetric matrix )
$$
Q_{\lambda }(x_1, \ldots , x_n)=
\begin{cases}Pf(Q_{\lambda_i\lambda_j })&\text{ if $m$ is even, }\cr
Q_{(\lambda, 0)}&\text{ if $m$ is odd.}\end{cases}\eqno(1.2.8)
$$

\ssbegin{1.2.3}{Lemma} Set
$\delta=\prod\limits_{1\leq i<j\leq n}\left(
\ds\frac{ e^{\frac12\varepsilon_{ij}}+ e^{\frac12\varepsilon_{ji}}}
{ e^{\frac12\varepsilon_{ij}}- e^{\frac12\varepsilon_{ji}}}\right)$, 
then

{\em i)}
$\delta^{-1}\Omega_3\delta=\sum\limits_1^n\pp_i^3-(\sum\limits_1^n\pp_i)^2. $

{\em ii)} The polynomials $Q_{\lambda}$ are common eigenfunctions of
the operators $\Omega_k$ for $k=1, 3, 5, \dots $.

{\em iii)} Let $P(x_1, \ldots , x_n)$ be polynomial symmetric with
respect to $ x_1, \ldots , x_n $ and such that after substitution
$x_i=t, \ x_j=-t$ it becomes independent of $t$. If $P$ is an
eigenfunction of all the operators $\Omega_k, \ k=1, 3, 5, \ldots$, then, 
up to a scalar multiple, $P$ coinsides with one of the $Q_{\lambda}$. 
\end{Lemma}

\ssec{1.3. Bispherical functions} Let $\fg$ be a finite dimensional
Lie superalgebra. Its enveloping algebra $U(\fg)$ possesses a
canonical antiautomorphism $t: u\mapsto \ ^tu$ which extends the {\it
principal antiautomorphism} of $\fg$ given by the formula $t(x)=-x$
for any $x\in \fg$ as follows:
$$
^t(uv)=(-1)^{p(u)p(v)}\ ^t(v)^t(u).
$$
\begin{Convention} For brevity, when it can not cause 
misunderstanding I sometimes write $(-1)^u$ 
instead of $(-1)^{p(u)}$. So $(-1)^{u+\bar 1}=(-1)^{p(u)+1}$ is not the same as 
$(-1)^{p(u+1)}$. \end{Convention}

The left and right coregular representations of $U(\fg)$ are defined
for any $l\in U(\fg)^*$ and $v, u \in U(\fg)$ by the formulas
$$
(L^*(u)l)(v)=(-1)^{ul}l(u^t\cdot v)\text{ and
$(R^*(u)l)(v)=(-1)^{u(l+v)}l(vu)$ }.
$$

Let $\fb_{1}$ and $\fb_{2}$ be subalgebras of $\fg$. A functional
$l\in U(\fg )^*$ is called {\it two-side invariant} if
$$
l(x_1u)=l(ux_2)=0\text{ for any $x_1\in\fb_1$, $x_2\in\fb_2$ and $u\in
U(\fg)$.}
$$
 
Let $V$ be a $\fg$-module containing a nonzero $\fb_{2}$-invariant
vector $v\in V$; suppose also that there exists a $\fb_{1}$-invariant
vector $v^*\in V^*$. The matrix coefficient $\Theta(v^*, v)\in
U(\fg)^*$ defined by the formula
$$
\Theta(v^*, v)(u)=(-1)^{uv}v^*(uv)\ \text{ for } \ u\in U(\fg)
$$
is a {\it bispherical function associated with the triple $(V, v^*,
v)$}.

Observe that if $z\in Z(\fg)$, then $L^*(z)l$ is two-side invariant if
so is $l$. Therefore, on the space of invariant functionals, every
$z\in Z(\fg)$ determines a linear operator $\Omega_{(z)}: l\mapsto
L^*(z)l$.

\ssec{1.4} In addition to the usual $\fg$-module structure on
$U(\fg)$, the extention of the adjoint action, consider the
$\fg$-module structure on $U(\fg)$ with respect to the following
action (\cite{Se4}, \cite{G})
$$
x*u=xu-(-1)^{x(u+\bar1)}ux.\eqno(1.4.1)
$$
It is easy to verify that for a finite dimensional $\fg$-module $V$
the functional $u\mapsto \tr_V(u)$ is invariant with respect to this
action.

In $\fg\oplus \fg$, consider two subalgebras:
$$
\fg_1=\{(x, (-1)^{x}x)\mid x\in \fg\}\text{ and }
\fg_2=\{(x, x)\mid x\in \fg\}.
$$
\begin{Lemma} The algebra of functionals on $U(\fg\oplus \fg)$
biinvariant with respect to $\fg_{1}$ and $\fg_{2}$ is isomorphic to
the algebra of functionals on $U(\fg)$ invariant with respect to the
action $(1.4.1)$.  \end{Lemma}

\ssbegin{1.4.1}{Convention} Let $I=I_{\bar 0}\cup I_{\bar 1}=\{1,
\ldots , n\}\cup\{\bar 1, \ldots , \bar n\}$ be the union of the
``even'' and ``odd'' indices.  Let $\dim V=(n|n)$ and $\{e_i\}_{i\in I}$
a basis of $V$ such that the parity of each vector of the basis is the
same as that of its index.   We assume that $\bar{\bar i}=i$ and define
the odd operator $\Pi\in \fgl(V)$ by setting
$$
\Pi(e_i)=(-1)^{p(i)}e_{\bar i}, \text{ for any } \ i\in I.
$$
Define the parity operator $P\in \fgl(V)$ by setting
$$
P(e_i)=(-1)^{p(i)}e_{i}, \text{ for any } \ i\in I.
$$

We denote the superspace of operators in $V$ by $\End(V)$, the
superspace of matrices (in the standard format with respect to the
basis in which all the even vectors come first) by $\Mat(V)$ and the
Lie superalgebra structure in these isomorphic superspaces is denoted
by $\fgl(V)$ or by $\fgl(\dim V)$.
\end{Convention}

\ssec{1.4.2. Queer superalgebras $\fq(V)$} Let
$$
\fq(n)=\{X\in \fgl(V)\mid [X, \Pi]=0\}.\eqno(1.4.2)
$$
It is easy to verify that $\fq(n)=\Span(e_{ij}, f_{ij}\mid i, j\in
I_{\bar 0})$, where 
$$
e_{ij}=e_i\otimes e_j+e_{\bar i}\otimes
e_{\bar \jmath}, \quad f_{ij}=e_i\otimes e^*_{\bar \jmath}+e_{\bar
i}\otimes e^*_j, \text{ for any }\ i, j\in I_{\bar 0} 
$$ 
and where $\{e_i^*\}_{i\in I}$ is the left dual of the basis
$\{e_i\}_{i\in I}$; we have also identified $\End(V)$ with $V\otimes
V^*$.

Set $\fg=\fq(n)\oplus \fq(n)$ and select $\fg_1$ and $\fg_2$ as at the
beginning of sec. 1.4. Let $\sigma: \fq(n)\tto\fg$ be the embedding
into the first summand, $\fh_{\bar 0}=\Span(e_{ii}\mid i\in I_{\bar
0})$ be the even part of Cartan subalgebra of $\fq(n)$. Let us define
inductively the following elements of $U(\fq(n))$:
$$
e_{ij}^{(1)}=e_{ij}, \ f_{ij}^{(1)}=f_{ij}, 
$$
$$
e_{ij}^{(p)}=\sum\limits_{l=1}^n e_{il}e_{lj}^{(p-1)}+(-1)^{k-1}
\sum\limits_{l=1}^n
f_{il}f_{lj}^{(p-1)}, \eqno (1.4.3)
$$
$$
f_{ij}^{(p)}=\sum\limits_{l=1}^n e_{il}f_{lj}^{(p-1)}+(-1)^{k-1}
\sum\limits_{l=1}^n
f_{il}e_{lj}^{(p-1)}.
$$

The following relations are subject to straightforward verification:
$$
\renewcommand{\arraystretch}{1.4}
\begin{array}{rcl}
[e_{ij}, e_{kl}^{(p)}]&=&\delta_{jk}e_{il}^{(p)}-\delta_{il}e_{kj}^{(p)}, 
\\[10pt]
[e_{ij}, f_{kl}^{(p)}]&=&\delta_{jk}f_{il}^{(p)}-\delta_{il}f_{kj}^{(p)}, 
\\[10pt] [f_{ij}, e_{kl}^{(p)}]&=&(-1)^{p+1}\delta_{jk}f_{il}^{(p)}
-\delta_{il}f_{kj}^{(p)}, \\[10pt]
[f_{ij}, f_{kl}^{(p)}]&=&(-1)^{p+1}\delta_{jk}e_{il}^{(p)}+
\delta_{il}e_{kj}^{(p)}.
\end{array}
\eqno(1.4.4)
$$
As is not difficult to verify, see \cite{Se3}, the elements
$z_k=\sum\limits_{i=1}^ne_{ii}^{(k)}, \ k=1, 3, 5, \ldots$, are central
ones if we embed $\fq(n)$ into $\fg=\fq(n)\oplus \fq(n)$ as the first
summand, i.e., strictly speaking, I mean not $z_k$ but 
$\sigma(z_k)\in U(\fg)$.

Let $\{\varepsilon_i\}_{i\in I_{\bar 0}}$ be the basis of $\fh_{\bar
0}$ dual to $\{e_{i}\}_{i\in I_{\bar 0}}$.  Let $\lambda= (\lambda_1,
\ldots , \lambda_m)$ be a strict partition of $k$, i.e.,
$\lambda_1>\lambda_2>\ldots>\lambda_m>0$.  Let $V^{\lambda}$ be an
irreducible submodule of $V^{\otimes k}$ corresponding to $\lambda $,
see \cite{Se1}.  Then the $\fg$-module
$w^{\lambda}=V^{\lambda}\otimes(V^{\lambda })$ is irreducible and
contains a unique even $\fg_2$-invariant vector $\omega_{\lambda}$
corresponding to the identity operator under identification
$V^{\lambda}\otimes(V^{\lambda})^*=\End(V^{\lambda})$.  The dual
module $(V^{\lambda})^*\otimes V^{\lambda}$ is also irreducible and
contains an even $\fg_1$-invariant vector $\omega^*_{\lambda}$
corresponding to the parity operator $P$ under identification
$(V^{\lambda})^*\otimes V^{\lambda}=\End(V^{\lambda})$.

\ssbegin{1.5}{Theorem} Let $\fq(n)$ be embedded into $\fg=\fq(n)\oplus
\fq(n)$ as the first summand. Then

{\em i)} Each two-side invariant with respect to $\fg_1$ and $\fg_2$
functional on $\fg$ is uniquely determined by its restriction onto
$S(\fh_{\bar 0})\subset \fq(n)$.

{\em ii)} Every $z\in Z(\fq(n))$ uniquely determines a differential
operator $\Omega_{(z)}$ on the space of restrictions of left-invariant
functionals $S(\fh_{\bar 0})^{*inv}$.

{\em iii)} $\Omega_{(z_k)}=\Omega_k$ for $z=z_{k}$.

{\em iv)} Let $\varphi_{\lambda}=\Theta(\omega^*_{\lambda}, 
\omega_{\lambda})$. Then the restriction of $\varphi_{\lambda}$ onto
$S(\fh_{\bar 0})$ coincides up to a scalar multiple with
$Q_{\lambda}(e^{\varepsilon_1}, \ldots , e^{\varepsilon_n})$
\end{Theorem}

\ssec{1.6} Let $\fg=\fgl(\dim V)$, where $\dim V=(n|n)$. Let
$I=I_{\bar 0}\cup I_{\bar 1}$, where $ I_{\bar 0}=\{1, \ldots , n\}$
and $ I_{\bar 1}=\{\bar 1, \ldots , \bar n\}$. Let $\{e_i\}$ be a
basis of $V$ such that the parity of each vector coincides with that
of its index. Denote by $\fpe_{1}(n)$ and $\fpe_{2}(n)$ the Lie
subsuperalgebras in $\fg$ preserving the respective tensors:
$$
\sum\limits_{i\in I}e^*_i\otimes e^*_{\bar i}, 
\quad
 \sum\limits_{i\in I}(-1)^{i}e^*_i\otimes
e^*_{\bar i}.
$$
Let $\psi_1$ and $\psi_2$ be involutive antiautomorphisms of $\fg$
that single out $\fpe_{1}(n)$ and $\fpe_{2}(n)$, respectively, i.e., 
$$
\fpe_{i}=\{X\in\fg l(v)\mid \psi_i(X)=-X\}\ \text{ for } i=1, 2.
$$
It is easy to see that the restrictions of $\psi_1$ and $\psi_2$ onto
$\fg_{\bar 0}$ coincide. Set
$$
\fh^{+}=\{x\in \fh\mid \psi_1(x)=\psi_2(x)=X\}, 
$$
where $\fh$ is Cartan subalgebra in $\fg$.

Let $\lambda$ be a partition and $V^{\lambda}$ the corresponding
$\fg$-submodule in the tensor algebra of the identity representation, 
see \cite{Se1}.

It follows from \cite{Se4} that, {(in Frobenius's notations \cite{Ma})}
$$
V^{\lambda}\text{ contains a $\fpe_{1}$-invariant vector if 
$\lambda=(\al_1, \ldots , \al_p|\al_1+1, \ldots , \al_p+1)$}
$$
and
$$
V^{\lambda}\text{ contains a $\fpe_{2}$-invariant vactor if 
$\lambda=(\al_1+1, \ldots , \al_p+1|\al_1, \ldots , \al_p)$.}
$$
Similarly, 
$$
(V^{\lambda})^*\text{ contains a $\fpe_{1}$-invariant vactor if 
$\lambda=(\al_1+1, \ldots , \al_p+1|\al_1, \ldots , \al_p)$}
$$
and 
$$
(V^{\lambda})^*\text{ contains a $\fpe_{2}$-invariant vactor if 
$\lambda=(\al_1, \ldots , \al_p |\al_1+1, \ldots , \al_p+1)$.}
$$
Now, let $\lambda=(\al_1+1, \ldots , \al_p+1|\al_1, \ldots , \al_p)$, let
$v_{\lambda}^*\in(V^{\lambda})^*$ be a $\fpe_{1}$-invariant vector, 
$v_{\lambda}\in(V^{\lambda})$ be a $\fpe_{2}$-invariant vector, and
let $\varphi_{\lambda}=\Theta(v_{\lambda}^*, v_{\lambda})$ be the
corresponding bispherical function. In $\fh^+$, select the basis
$\{e^+_{ii}=\frac12(e_{ii}+\psi_1(e_{ii}))\mid i\in I_{\bar 0}\}$ and let
$\varepsilon_{i}$ be the left dual vectors.

Let $e_{ij}$ be the basis of matrix units in $\fgl(V)$. Set
$$
\eij^{(1)}=\eij, \quad \eij^{(k)}=\mathop{\sum}\limits_{i\in 
I}(-1)^pe_{ip}e_{pj}^{(k-1)}.
$$
It is easy to verify that 
$$
[e_{ij}, 
e_{pq}^{(k)}]=\delta_{jp}e_{ip}^{(k)}-(-1)^{(i+j)(p+q)}\delta_{iq}e_{pj}^{(k)}.
$$
This easily implies that $z_{k}=\mathop{\sum}\limits_{i\in 
I}e_{ii}^{(k)}\in Z(\fg)$.

\ssbegin{1.7}{Theorem} {\em i)} Each two-sided invariant with respect 
to $\fpe_1$ and $\fpe_2$ functional on $U(\fg)$ is uniquely determined by its 
restriction onto $S(\fh^+)$.

{\em ii)} Every $z\in Z(\fg)$ uniquely determines a differential
operator $\Omega_{(z)}$ on the space $S(\fh^+)^*$ of restrictions of
invariant functionals.

{\em iii)} If $z=z_k$, then $\Omega_{(z_k)}=\tilde {\Omega}_k$.

{\em iv)} The functional $\varphi_{\lambda}$ coincides, up to a scalar
multiple, with $Q_{\lambda}(e^{\varepsilon_1}, \ldots , 
e^{\varepsilon_k})$.
\end{Theorem}

\section*{\S 2. The algebra dual to the enveloping algebra}

In this section we follow Dixmier's book \cite{Dix} applied {\it
mutatis mutandis} to Lie superalgebras.

\ssec{2.1} Let $\fg$ be a Lie superalgebra. We endow $U(\fg)^*$ with
a coalgebra structure by setting
$$
c:\fg \tto U(\fg)\otimes U(\fg), \quad c(x)=x\otimes 1+1\otimes x \ 
\text{ for any }x\in \fg, 
$$
so that $^tc(x)=c(^tx)$ where the first $t$ is the principal
antiautomorphism of the superalgebra $U(\fg)\otimes U(\fg)=
U(\fg\oplus\fg)$, see sec. 1.3.

\ssbegin{2.2}{Lemma} Let $\dim \fg=(n|m)$. Then $ U(\fg)^*$ is
isomorphic to the supercommutative superalgebra of formal power series
in $n$ even and $m$ odd indeterminates. \end{Lemma}

\begin{proof} Let $ \fg_{\bar 0}=\Span(e_1, \ldots , e_n)$, 
$\fg_{\bar 1}=\Span(e_{\bar 1}, \ldots , e_{{\overline m}})$, $I_{\bar
0}=\{1, \ldots , n\}$, $ I_{\bar 1}=\{\bar 1, \ldots , \overline m\}$.
 
Denote: 
$$
M=\{\nu=(\nu_1, \ldots , \nu_n, \nu_{\bar i}, \ldots , \nu_{{\overline 
m}})\mid \nu_i\in \Zee_{\ge 0}\text{ if $i\in I_{\bar 0}$ and
$\nu_i\in\{0, 1\}$ if $i\in I_{\bar 1}$} \}.
$$ 
For $\nu\in M$, set 
$$
e_{\nu}=\ds\frac{e^{v_1}_1}{\nu_1!}\cdots \ds\frac{e^{v_n}_n}{\nu_n!}
\ds\frac{e^{v_{\bar 1}}_{\bar 1}}{(\nu_{\bar 1})!}\cdots
\ds\frac{e^{v_{{\overline m}}}_{{\overline m}}}{(\nu_{{\overline m}})!}
$$ 
and let $t_1, \ldots , t_n, t_{\bar 1}, \ldots , t_{{\overline m}}$ be the set of
even and odd supercommuting indeterminates. The correspondence
$$
U(\fg)^*\ni L\mapsto \sum\limits_{\nu\in M}L(e_{\nu})t^{\nu_{{\overline m}}}_{{\overline m}}
\ldots t^{\nu_{\bar 1}}_{\bar 1} t^{\nu_{n}}_{n}\ldots t_1^{\nu_1}
$$
determines the homomorphism desired. \end{proof}

\ssec{2.3. Left and right coregular representations} For any $u, v\in
U(\fg)$ and $L\in U(\fg)^*$ set
$$
(L^*(u)L)(v)=(-1)^{uL}L(^tuv)\text{ and }
(R^*(u)L)(v)=(-1)^{u(L+v)}L(vu).
$$
The following statements are easy to verify:

i) $u\tto L^*(u)$ is a representation of $U(\fg)$ in $U(\fg)^*$ (we
call it the {\it left regular} one);

ii) $u\tto R^*(u)$ is a representation of $U(\fg)$ in $U(\fg)^*$ (we
call it the {\it right regular} one);

iii) If $x\in \fg$, then $L^*(x)$ and $R^*(x)$ are
superdifferentiations of the algebra $U(\fg)^*$.

\noindent Observe also that algebra $U(\fg)^*$ possesses an
automorphism
$$
L\mapsto L^T: L^T(u)=L(^tu)\ \text{ for any }\ u\in U(\fg), \ L\in
U(\fg)^*.
$$

\ssec{2.4. Matrix coefficients} Let $V$ be a $\fg$-module, $V^*$ the
dual module, let $v\in V$ and $v^*\in V^*$. Let $\pi: \fg\tto\fgl(V)$
be the corresponding representation. Denote by $\Theta^{\pi}(v^*, v)$
the linear form on $U(\fg)^*$:
$$
\Theta^{\pi}(v^*, v)(u)=(-1)^{uv}v^*(\pi(u)v). \eqno(2.4.1)
$$
Finally, denote by $C(\pi)$ (or $C(V)$) the subalgebra in $U(\fg)^*$
generated by $\Theta^{\pi}(v^*, v)$ for all $v^*\in V^*$ and $v\in V$.

\ssbegin{2.4.1}{Lemma} 
$$
\Theta^{\pi_1\ot \pi_2}(v_1^*\ot v_2^*, v_1\ot
v_2)=(-1)^{v_1 v_2^*}\Theta^{\pi_1}(v_1^*, v_1)\ot
\Theta^{\pi_2}(v_2^*, v_2)
$$
\end{Lemma}

\begin{proof} We have 
$$
\renewcommand{\arraystretch}{1.4}
\begin{array}{rl}
&\Theta^{\pi_1\otimes \pi_2}(v_1^* \otimes v_2^*, v_1\otimes v_2)
(u_1\otimes u_2) \\
[10pt] =&(-1)^{(u_1+u_2)(v_1+v_2)} v_1^* \otimes
v_2^*(\pi_1(u_1)\otimes\pi_2(u_2)(v_1\otimes v_2)) \\
[10pt] =&
(-1)^{(u_1+u_2)(v_1+v_2)+v_1u_2+v_2^*(u_1+v_1)}
v_1^*(\pi_1(u_1)v_1)v_2^*(\pi_2(u_2)v_2) \\
[10pt]
=&\Theta^{\pi_1}(v_1^*, v_1)\otimes
\Theta^{\pi_2}(v_2^*, v_2)(u_1\otimes u_2) \\
[10pt] =&
(-1)^{u_1(v_2+v_2^*)} \Theta^{\pi_1}(v_1^*, v_1)(u_1)
\Theta^{\pi_2}(v_2^*, v_2)(u_2) \\
[10pt]
=&(-1)^{u_1(v_2+v_2^*)+v_1u_1+v_2u_2}
v_1^*(\pi_1(u_1)v_1)v_2^*(\pi_2(u_2)v_2).
\end{array}
$$
\end{proof} 

\ssbegin{2.4.2}{Lemma} The map $V^*\otimes V\tto U(\fg)$ given by $v^*\otimes
v\mapsto\Theta^{\pi}(v^*, v)$ is a $\fg\oplus\fg$-module homomorphism;
here we consider $U(\fg)^*$ as a $\fg\oplus\fg$-module with respect to
the left and right coregular representations. If $V$ is irreducible, 
the above map is an isomorphism. \end{Lemma}

\begin{proof} Clearly, there exists a linear map $\varphi:V^*\otimes
V\tto U(\fg)^*$ such that $\varphi(v^*\otimes v)=\Theta(v^*, v)$. Let
$x\in\fg$. Then
$$
\renewcommand{\arraystretch}{1.4}
\begin{array}{l}
\varphi\left((x\ot 1)(\vz\ot v)\right)(u)=\vh(x\vz\ot
v)(u)=\\
\ttt(x\vz, v)(u)=(-1)^{vu}(x\vz)(u)=(-1)^{vu+\bar1+x\vz}\vz(xu).
\end{array}
$$
On the other hand, 
$$
\renewcommand{\arraystretch}{1.4}
\begin{array}{rcl}
L^*(x)\vh(\vz, v)(u)
&=&
(-1)^{x(v+\vz)+\bar1}\ttt(\vz, v)(xu)
\\[10pt]
&=&(-1)^{x(v+\vz)+\bar1+(x+u)v}\vz(xuv)
\\[10pt]
&=&(-1)^{uv+\bar1+x\vz}\vz(xu).
\end{array}
$$
Further, 
$$
\vh((1\ot x)(\vz\ot v))(u)=(-1)^{x\vz}\vh(\vz\ot xv)(u)=
(-1)^{x\vz}\ttt(\vz, xv)(u)=
(-1)^{x\vz+u(x+v)}\vz(uxv).
$$
On the other hand, 
$$
\renewcommand{\arraystretch}{1.4}
\begin{array}{l}
 R^*(x)\vh(\vz\ot v)(u)=
(-1)^{x(\vz+v+u)}\vh(\vz\ot v)(ux)=
(-1)^{x(\vz+v+u)}\ttt(\vz\ot v)(ux)=\\
(-1)^{x(\vz+v+u)+(u+x)v}\vz(uxv)=
(-1)^{x\vz+u(x+v)}\vz(uxv).\end{array}
$$
This proves the first statement. The second one is obvious.
\end{proof} 

\ssbegin{2.4.3}{Lemma} Let $\rho$ be the tensor product of 
representations $\pi_1$ and $\pi_2$. Then

{\em i)} $\ttt^{\rho}(\vv\ot\vvv, v_1\ot v_2)=(-1)^{v_1\vvv}\ttt^{\pi_1}
(\vv, v_1)\ttt^{\pi_2}(\vvv, v_2)$.

{\em ii)} $C(\rho)=C(\pi_1)C(\pi_2)$.

{\em iii)} If $\pi$ is finite dimensional, then
$(\ttt^{\pi}(\vz, v))^T=(-1)^{v\vz}\ttt^{\pi^*}(v, \vz)$. \end{Lemma}

\begin{proof} i)
$$
\renewcommand{\arraystretch}{1.4}
\begin{array}{rcl}
\ds(-1)^{v_1\vvv}(\ttt^{\pi_1}(\vv, v_1)\ttt^{\pi_2}(\vvv, v_2))(u)
&=&
(-1)^{v_1\vvv}(\ttt^{\pi_1}(\vv, v_1)\ot\ttt^{\pi_2}(\vvv, v_2))(C(u))
\\[10pt]
&=&
(-1)\ttt^{\pi_1\ot\pi_2}(\vv\ot\vvv, v_1\ot v_2)(C(u))
\\[10pt]
&=&
(-1)^{(v_1+v_2)u}(\vv\ot\vvv)(\pi_1\ot\pi_2\cdot C(u)(v_1\ot v_2)
\\[10pt]
&=&
(-1)^{(v_1+v_2)u}(\vv\ot\vvv)(\rho(u)(v_1\ot v_2))
\\[10pt]
&=&
\ttt^{\rho}(\vv\ot\vvv, v_1\ot v_2)(u).
\end{array}
$$
ii) follows from i)

iii) 
$$
\renewcommand{\arraystretch}{1.4}
\begin{array}{rcl}
(\ttt^{\pi}(\vz, v))^T(u)
&=&\ttt^{\pi}(\vz, v)(^tu)
\\[10pt]
&=&(-1)^{uv}\vz(\pi(^tu)(v)).
\end{array}
$$
On the other hand, 
$$
\renewcommand{\arraystretch}{1.4}
\begin{array}{rcl}
\ttt^{\pi^*}(v, \vz)(u)
&=&
(-1)^{\vz u}v(\pi^*(u)\vz)
\\[10pt]
&=&
(-1)^{\vz u+v(u+\vz)}(\pi^*(u)\vz)(v
)
\\[10pt]
& =&
(-1)^{\vz u+v(u+\vz)+u\vz}\vz(\pi(^tu)v).
\end{array}
$$
\end{proof} 

\ssbegin{2.4.4}{Lemma} Let $V$ be a $\fg$-module. Consider the map 
$$
\renewcommand{\arraystretch}{1.4}
\begin{array}{l}
\vh:V^*\ot V\ot V^* \ot V\tto U(\fg)^*\\
\vh(\vv\ot v_1\ot\vvv\ot v_2)=
(-1)^{\vv\vvv+\vvv v_1+\vv v_1}(\ttt(\vvv, v_1))^T\cdot \ttt(\vv, v_2).\end{array}
$$
If we consider $V^*\ot V\ot V^* \ot V=(V^*\ot V)\ot (V^* \ot V)$ as a
$\fg\oplus\fg$-module in such a way that the first two factors is one
$\fg$-module and the last two ones is the other module, then $\varphi$
is a $\fg\oplus\fg$-module homomorphism. (Recall that we consider
$U(\fg)^*$ as a $\fg\oplus\fg$-module with respect to the simultaneous
left and right coregular representations.) \end{Lemma}

\begin{proof} Let $x\in \fg$. Then
$$
\renewcommand{\arraystretch}{1.4}
\begin{array}{l}
\varphi((x\otimes 1)(v_1^*\otimes v_1\otimes v_2^*\otimes v_2)) =
\varphi((xv_1^*\otimes v_1+(-1)^{xv_1^*}v_1^*\otimes xv_1) \otimes
v_2^*\otimes v_2)) \\
[10pt] =
(-1)^{v_2^*(x+v_1^*)+v_2^*v_1+(x+v_1^*)v_1}
(\Theta(v_2^*, v_1))^T\Theta(xv_1^*, v_2) +\\
[10pt] 
(-1)^{v_1^*v_2^*+v_2^*(x+v_1)+v_1^*(x+v_1)+xv_1^*}(\Theta(v_2^*, xv_1))^T
\Theta(v_1^*, v_2) \\
[10pt]
=(-1)^{v_1^*v_2^*+v_2^*v_1+v_1^*v_1}L^*(x)
(\Theta(v_2^*, v_1))^T\Theta(v_1^*, v_2) +\\
[10pt] 
(-1)^{v_1^*v_2^*+v_2^*v_1+v_1^*v_1+x(v_2^*+v_1)}(\Theta(v_2^*, v_1))^T
L^*(x)\Theta(xv_1^*, v_2) \\
[10pt] =
(-1)^{v_1^*v_2^*+v_2^*v_1+v_1^*v_1}L^*(x)[(\Theta(v_2^*, v_1))^T
\Theta(v_1^*, v_2)] \\[10pt] = L^*(x)\varphi(v_1^*\otimes v_1\otimes
v_2^*\otimes v_2).
\end{array}
$$

The identity 
$$
\varphi((1\otimes x)(v_1^*\otimes v_1\otimes v_2^*\otimes v_2))=
R^*(x)\varphi(v_1^*\otimes v_1\otimes v_2^*\otimes v_2)
$$
is similarly verified.
\end{proof} 

\ssbegin{2.4.5}{Lemma} Let $V$ be a finite dimensional $\fg$-module,
$\{v_i\}_{i\in I}$ its basis, $\{v_i^*\}_{i\in I}$ the left dual
basis of $V^*$.  Then
$$
\sum\limits_i(\ttt(\vz, \vi))^T\ttt(\vzi, v)=\vz(v)\cdot\vr, 
$$
where $v\in V$, $\vz\in V^*$ and $\vr\in U(\fg)^*$ is the counit.
\end{Lemma}

\begin{proof} The functional $\vr$ is uniquely, up ot a scalar
multiple, characterized by its invariance with respect to the right
coregular representation. Further, 
$$
\omega=\sum (-1)^{i}\vzi\ot\vi
$$
is an invariant of the $\fg$-module $V^*\ot V$. Hence, by Lemma
2.4.4, 
$$
\vh(\sum(-1)^{i}\vzi\ot\vi\ot\vz\ot v)
$$ 
is an invariant with respect to the right coregular representation on
$U(\fg)$, i.e., $\vh$ is an invariant functional on $U(\fg)$. Hence,
$\vh(\omega\ot\vz\ot v)=\al\vr$.

On the other hand, 
$$
\renewcommand{\arraystretch}{1.4}
\begin{array}{l}
\vh(\omega\ot\vz\ot v)=\sum(-1)^{i)}\vh(\vzi\ot\vi\ot\vz\ot v)=\\
\sum(\ttt(\vz, \vi))^T\cdot\ttt(\vzi, v).\end{array}
$$
Hence, 
$$
\sum\limits_i (\ttt(\vz, \vi))^T\cdot\ttt(\vzi, v)=\al\vr.
$$
To find $\al$, let us substitute $u=1$ into both parts of the
identity. We obtain:
$$
\al=\al\vr(1)=\sum((\ttt(\vz, \vi))^T\cdot\ttt(\vzi, v)(1)=
\sum \vz(\vi)\vzi(v)=\vz(v).
$$
\end{proof} 

Let $\fg=\fg l(V)$ and let $\fA$ be the subalgebra in $U(\fg)^*$
generated by $C(V)$ and $C(V^*)$. It is not difficult to verify that
$\fA$ is invariant with respect to the left and right coregular
representations.

\section*{\S 3. Bispherical functions and radial parts of Laplace
operators for the triple $(\fq(n), \fq(n)\oplus \fq(n), \fq(n))$}

\ssec{3.1} Let $\fg=\fq(n)$, let $Z(\fg)$ be the center of $U(\fg)$
and $\fh=\Span (e_{ii}, f_{jj}\mid i, j\in I_{\ev})$ the
Cartan subalgebra (for definition of $\eij$ and $\fij$ see sec. 
1.4.2). Let $\mu\in\fh^*_{\bar 0}$ and $M(\mu)$ the coresponding
Verma module with the highest weight with respect to the decomposition
$\fg=\fn_{-}\oplus\fh\oplus\fn_{+}$, where
$\fn_-=\Span(\eij, \fij\mid i>j)$ and
$\fn_+=\Span(\eij, \fij)\mid i<j)$, see \cite{Se4}, 
\cite{Pe}. Such module $M(\mu)$ has a central character, i.e., there
exists a homomorphism
$$
\chi_{\mu}:Z(\fg)\tto \Cee,\quad zv=\chi_{\mu}(z)v\text{ for any $v\in
M(\mu)$ and $z\in Z(\fg)$}.
$$
Let $\eps_{i}$ be the left dual vector to $e_{ii}$. Set
$R^+=\{\eps_{i}-\eps_{j}\mid i<j\}$ and
$$
B(\mu)=\{\nu\mid \nu=\lambda-\sum n_{\al}\al \text{ for any 
$n_{\al}\in \Zee_{\geq 0}$, where $\al\in R^+$, and
$\chi_{\nu}=\chi_{\mu}$}\}. 
$$

\ssbegin{3.2}{Lemma} Let $L(\mu)$ be an irreducible module with
highest weight $\mu$. Then
$$
\ch L(\mu)=\sum\limits_{\nu\in B(\mu)} C_{\mu\nu}\, \ch M(\nu)
\eqno(3.2.1)
$$ 
where for any fixed $\tau\in\fh^*_{\bar 0}$ only a finite number of
summands in the right hand side contains $\tau$ as a weight. \end{Lemma}

\begin{proof} Let $\al_1, \ldots , \al_{n-1}$ be a base (system of simple
roots) in the root system $R$. Set
$$
ht(\nu)=\sum\limits_1^{n-1}k_i \quad \text{for any }
\nu=\mu-\sum\limits_{i=1}^{n-1}k_i\al_i.
$$
By induction on $t$ we prove that 
$$
\ch L(\mu)=\sum\limits_{\nu\in B(\mu), \; ht(\nu)\le t}
C_{\mu\nu}\; \ch M(\nu)+\sum\limits_{i=1}^{m(t)}\ch V_i
\eqno(3.2.2)
$$
for any $\tau\in \Supp\ \ch V_i$, $ht(\tau)>t$, 
$\chi_{V_i}=\chi_{\mu}$ and $\Supp\ \ch V_i\subset B(\mu)$.

Let $t=0$. Then we obtain an exact sequence
$$
0\tto N\tto M(\mu)\tto L(\mu)\tto 0
$$
with $N\subset \mathop{\oplus}\limits_{\nu\ne \mu} M(\mu)_{\nu}$. Therefore, 
$ht(\nu)>0$ and $\ch L(\mu)=\ch M(\mu)-\ch(N)$; moreover, $\Supp\ \ch
N\subset B(\mu)$.

Let (3.2.2) hold. Consider the set $ht(\nu)$, where $\nu\in \Supp\
\ch V_i$ and let $k$ be the least element of this set. Clearly, 
$k>t$. There exists then a finite number of vectors $v_1^i, \ldots , 
v_l^i\in V_i$ with weights $\mu_{1i}, \dots, \mu_{li}$ such that
$ht(\mu_{1i})=ht(\mu_{2i})\dots =k$. For all the other weights $\nu$
of $V_i$ we have $ht(\nu)>k\ge t+1$. The vectors $v_1^i, \ldots , v_l^i$
are, obviously, the highest weight ones, so we have an exact sequence
$$
0\tto N\tto \mathop{\oplus}\limits_{j=1}^l
 M(\mu_{ji})\tto V_i\tto K\tto 0.
\eqno (3.2.3)
$$
This implies:
$$
\ch V_i=\sum\limits_{j=1}^l \ch M(\mu_{ji})+\ch K-\ch N
$$ 
and if $\nu\in\Supp \ \ch K\cup\Supp\ \ch N$, then $ht(\nu)>t+1$. 
This proves (3.2.2).

Let $\nu\in\fh^*_{\bar 0}$ be any weight. Set $t=ht(\nu)$; apply
(3.2.2) to see that $\nu\notin \cup\Supp\ \ch V_i$. This proves
(3.2.1).
\end{proof}

\ssbegin{3.3}{Lemma} Let $z_3=\sum\eii^{(3)}$ (see $1.4.2$) and let
$v_{\mu}$ be the highest weight vector of the $\fg$-module $M(\mu)$. 
Then
$$
z_3v_{\mu}=\left(\sum \mu_i^3-(\sum\mu_i)^2\right)v_{\mu}.
$$ \end{Lemma}

\begin{proof} First, observe that
$$
\renewcommand{\arraystretch}{1.4}
\begin{array}{l}
 \eii^{(3)}=\sum\limits_{j=1}^n \eij\eji^{(2)}+ \sum\limits_{j=1}^n 
 \fij\fji^{(2)}= \sum\limits_{i<j} \eij\eji^{(2)}+ \sum\limits_{i<j} 
 \fij\fji^{(2)}\\
+\sum\limits_{i>j} \eij\eji^{(2)}+ \sum\limits_{i>j} \fij\fji^{(2)}
+\sum \eii\eii^{(2)}+ \sum \fii\fii^{(2)}\end{array}
$$
But, as is easy to verify, $\eij^{(2)}v_{\mu}=\fij^{(2)}v_{\mu}=0$. 
Therefore, 
$$
\renewcommand{\arraystretch}{1.4}
\begin{array}{l}
 z_3v_{\mu}=\left(\sum\eii\eii^{(2)}+\sum\fii\fii^{(2)}\right)v_{\mu}+
\left(\sum\limits_{i<j}\eij\eji^{(2)}+\sum\fij\fji^{(2)}\right)v_{\mu}\\
=\left(\sum\eii\eii^{(2)}+\sum\fii\fii^{(2)}+ \sum\limits_{i<j} \eii^{(2)}-
\ejj^{(2)}+ \sum\limits_{i<j} \ejj^{(2)}-
\eii^{(2)}\right) v_{\mu}\\
=\left(\sum\eii\eii^{(2)}+ \sum\fii\fii^{(2)}\right) v_{\mu}.\end{array}
$$
Further, it is easy to verify that
$$
\eii^{(2)} v_{\mu}=\left(\mu_i^2-\mu_i-2\sum\limits_{k>i}\mu_k\right)
 v_{\mu}
\text{ and } \fii^{(2)} v_{\mu}=0.
$$ 
Hence, $z_3v_{\mu}=(\sum\mu_i^3-(\sum\mu_i )^2)v_{\mu}$. \end{proof} 

\ssbegin{3.4}{Corollary} Let $\delta=\prod\limits_{\al\in
R^+}\ds\frac{e^{\al/2}+ e^{-\al/2}}{ e^{\al/2}- e^{-\al/2}}$ (see
$1.4.2$) and let $\Omega_3^*=\sum\pp_i^3-(\sum\pp_i)^2$, where
$\pp_ie^l=l(e_i)e^l$. Then $\Omega_3^*\delta^{-1}=0$ and $(\sum\pp_i^3)\delta^{-1}=0$.
\end{Corollary}

\begin{proof} Let us prove a more general statement; namely, let
$\vh_{\mu}$ be the character of an irreducible $\fg$-module with
highest weight $\mu$. Then 
$$
\delta^{-1}\vh_{\mu}\ \text{{\sl is an eigenfunction of $\Omega_3^*$
with eigenvalue $\sum\mu_i^3-(\sum\mu_i)^3$.}}
$$

Indeed, by Lemma 3.2 
$$
\vh_{\mu}=\sum\limits_{\nu\in B(\mu)}C_{\mu\nu} \ch M(\nu).
$$ 
By multiplying both parts of the inequality by $\delta^{-1}$ we obtain
$$
\delta^{-1}\vh_{\mu}=\sum\limits_{\nu\in B(\mu)}\tilde
C_{\mu\nu}e^{\nu}.\eqno(3.4.1)
$$ 
If $z\in Z(\fg)$, then $zv_{\mu}=\chi_{\mu}(z)v_{\mu}$. Since $\nu\in
B(\mu)$, it follows that $\chi_{\mu}(z)=\chi_{\nu}(z)$. But by Lemma
3.3
$$
\chi_{\mu}(z_3)=\sum\mu_i^3-\left(\sum\mu_i\right)^2
=\chi_{\nu}(z_3)=\sum\nu_i^3-\left(\sum\nu_i\right)^2.
$$
So for $\Omega_3^*=\sum\pp_i^3-(\sum\pp_i)^2$, we have
$$
\renewcommand{\arraystretch}{1.4}
\begin{array}{l}
 \Omega_3^*\delta^{-1}\vh_{\mu}=\sum\limits_{\nu\in B(\mu)}
\Omega_3^*(\tilde C_{\mu\nu}e^{\nu})=
\sum\limits_{\nu\in B(\mu)} \tilde C_{\mu\nu}\left(\sum \vi^3-\left(
\sum \vi\right)^2\right)e^{\nu}\\
=\left(\sum\mu_i^3-\left(\sum\mu_i\right)^2\right)\delta^{-1}\vh_{\mu}.
\end{array}
$$
In particular, applying this statement to the trivial module we
obtain $\Omega_3^*\delta^{-1}=0$, implying
$(\sum\pp_i^3)*\delta^{-1}=0$.
\end{proof}

\ssec{3.5.  Proof of heading i) of Lemma 1.2.3} Let
$R=\{\vr_i-\vr_j\mid i\ne j\}$ and $R^+=\{\vr_i-\vr_j\mid i< j\}$. 
Set $\de_{\al}^+=\eal+\ealm$ and $\de_{\al}^-=\eal-\ealm$ for any
$\al\in R$.  Then
$$
\pp_i(\de^+_{\al})=\frac12\al(\ei)\de_{\al}^-, \quad
\pp_i(\de^-_{\al})=\frac12\al(\ei)\de_{\al}^+.
$$ 
We also set
$$
\vh_i=\sum\limits_ {\al\in
R^+}\frac{\al(\ei)}{\de_{\al}^+\de_{\al}^-}, \qquad
\psi_i=\sum\limits_{\al\in R^+}\frac{\al^2(\ei)}{(\de_{\al}^+)^2},
\qquad \ttt_i=\sum\limits_{\{\al, \be\}\subset R^+}
\frac{\al(\ei)\be(\ei)} {\de_{\al}^+\de_{\al}^-
\de_{\be}^+\de_{\be}^-} \eqno(3.5.1)
$$ 
where the last sum runs over the two-element subsets of $R^+$.  It is
not difficult to verify that one can express the operator $\hat
{\Omega}_3=\Omega_3+(\sum\pp_i)^2$ in the form
$$
\hat {\Omega}_3=\sum\pp_i^3+6\sum\vh_i\pp_i^2-6\sum\psi_i\pp_i+24\sum\ttt_i\pp_i.
\eqno(3.5.2)
$$ 
It is easy to verify that $\delta^{-1}\pp_i\delta=\pp_i-2\vh_i$; hence, 
$$
\delta^{-1}\pp_i^2\delta =\pp_i^2-4\vh_i\pp_i+4\vh_i^2-2\pp_i(\vh_i)
\eqno(3.5.3)
$$ 
$$
\delta^{-1}\pp_i^3\delta
=\pp_i^3-6\vh_i\pp_i^2+3(4\vh_i^2-2\pp_i(\vh_i))\pp_i-
8\vh_i^3-2\pp_i^2(\vh_i)+12\vh_i\pp_i(\vh_i)\eqno(3.5.4)
$$
Therefore, 
$$
\renewcommand{\arraystretch}{1.4}
\begin{array}{l}
 \delta^{-1}\hat {\Omega}\delta
=\sum_i[\pp_i^3+(24\ttt_i-6\psi_i-12\vh_i^2-6\pp_i(\vh_i))\pp_i\\
+(16\vh_i^3-2\pp_i^2(\vh_i)-48\ttt_i\vh_i+12\psi_i\vh_i)].\end{array}
$$
Direct calculations show that 
$$
24\ttt_i-6\psi_i-12\vh_i^2-6\pp_i(\vh_i)=0.
$$ 
Hence, 
$$
\delta^{-1}\hat {\Omega}\delta =\sum\pp_i^3+f.
$$ 
But, $\hat {\Omega}(1)=0$, so 
$$
(\delta^{-1}\hat {\Omega}\delta )\cdot(\delta^{-1})=
\left(\sum\pp_i^3\right)( \delta^{-1})+f\delta^{-1}=0.
$$ 
But due to Corollary 3.4, $(\sum\pp_i^3)( \delta^{-1})=0$; hence, 
$f\delta^{-1}=0$ and $f=0$. \qed

\ssbegin{3.6}{Lemma} Let $\fg$, $\fg_1$, $\fg_2$ be selected as in
sec. $1.4$, let $I$ be the left ideal in $U(\fg\oplus \fg)$ generated
by $\fg_2$ and $M= U(\fg\oplus \fg)/I$. Let $\sigma:\fg\tto\fg\oplus
\fg$ be the embedding into the first summand, i.e., $\sigma(x)=(x,
0)$. Let $\tilde {\sigma}: U(\fg)\tto M$ be the map induced by the
homomorphism $U(\fg)\tto U(\fg\oplus \fg)$ that extends $\sigma$ and
$\rho(x)=(x, (-1)^xx)$ an isomorphism of $\fg$ with $\fg_1$. Then
$\tilde {\sigma}(x*u)=\rho(x) \tilde {\sigma}(u)$. \end{Lemma}

\begin{proof} 
$$
\renewcommand{\arraystretch}{1.4}
\begin{array}{rcl}
 \tilde {\sigma}(x*u)&=& \tilde {\sigma}(xu-(-1)^{x(u+\bar1)}ux)=xu\ot
 1-(-1)^{x(u+\bar1)}ux\ot 1 \\[10pt] &=& xu\ot
 1-(-1)^{x(u+\bar1)}ux\ot 1-\rho(x) \tilde {\sigma}(u)+\rho(x) \tilde
 {\sigma}(u) \\[10pt] &=&\rho(x) \tilde {\sigma}(u)+xu\ot
 1-(-1)^{x(u+\bar1)}ux\ot 1 \\[10pt] &&-(x\ot 1+(-1)^x1\ot x)(u\ot 1)
 \\[10pt] &=& \rho(x) \tilde {\sigma}(u)-(-1)^{x(u+\bar1)}(ux\ot
 1+u\ot x) \\[10pt] &=& \rho(x) \tilde
 {\sigma}(u)-(-1)^{x(u+\bar1)}(u\ot 1)(x\ot 1+1\ot x) \\[10pt]
 &\equiv_{\pmod I}& \rho(x) \tilde {\sigma}(u).
\end{array}
$$
\end{proof}

\ssec{3.6.1. Statement of Lemma 1.4  is true} This is a direct 
corollary of Lemma 3.6.

\ssec{3.7. Proof of heading i) of Theorem 1.5} Due to \cite{Se4} and 
\cite{G} we have an isomorphism of $\fg$-modules with respect to the 
action (1.4.1):
$$
U(\fg)\simeq\Ind_{\fg_{\bar 0}}^{\fg}(U(\fg_{\bar 0})).
$$ 
Therefore, there exists a bijection between the space of
$\fg$-invariant with respect to the action (1.4.1) functionals on
$U(\fg)$ and the space of $\fg_{\ev}$-invariant functionals on
$U(\fg_{\ev})$. Moreover, any $\fg$-invariant functional is uniquely
determined by its restriction onto $U(\fg_{\ev})$. On the other hand, 
every $\fg_{\ev}$-invariant functional is uniquely determined by its
restriction onto $U(\fh_{\ev})=S(\fh_{\ev})$. \qed

\ssec{3.8} Let $l\in (U(\fg)^*)^{\fg}$ and $\vh_l$ the generating
function of its restriction onto $S(\fh_{\ev})$, i.e., 
$$
\vh_l (t_1, \ldots , t_n)=\sum\frac{l(e^{\nu_1}_{11}\ldots e^{\nu_n}_{nn})}
{(\nu_1)!\ldots(\nu_n)!} t^{\nu_1}_{1}\ldots t^{\nu_n}_{n}.\eqno(3.8.1)
$$ 
On $S(\fh_{\ev})^*$, define the following operators by setting for any $f\in S(\fh_{\ev})$:
$$
(\pp_i^{(k)}l)(f)=l(\eii^{(k)}f), \quad
(\delta_i^{(k)}l)(f)=L(\fii^{(k)}f), \eqno(3.8.2)
$$ 
$$
(D_{ij}^{(k)}l)(f)=l(\eij\eji^{(k)}f), \quad
(\Delta_{ij}^{(k)}l)(f)=L(\fij\fji^{(k)}f).\eqno(3.8.3)
$$ 

\ssbegin{3.8.1}{Lemma} Let $\al=\vr_i-\vr_j$. Then

{\em i)}
$
D_{ij}^{(k)}=\ds\frac{e^{\al}}{ e^{\al}-1}(\pp_i^{(k)}- \pp_j^{(k)}).
$

{\em ii)}
$
\Delta_{ij}^{(k)}=\ds\frac{e^{\al}}{ e^{\al}+1}(\pp_j^{(k)}+(-1)^{k+1} \pp_i^{(k)}).
$

{\em iii)}
$l(f\cdot \fii\fii^{(k)})=l(f\eii^{(k)})$  for $k$ odd.

{\em iv)} $ l(f\cdot \fii\fii^{(k)})=0$
 for $k$ even.
\end{Lemma}

\begin{proof} i) and ii) are similarly proved. Consider i): 
$$
\renewcommand{\arraystretch}{1.4}
\begin{array}{l}
(D_{ij}^{(k)}l)(f)=l(f\eij\eji^{(k)})
=L(\eij f(h+\al(h))\eji^{(k)}=
l(f(h+\al(h))\eji^{(k)}\eij)
\\[10pt]
=l(f(h+\al(h))\eij\eji^{(k)})-
l(f(h+\al(h))[\eij, \eji^{(k)}])
\\[10pt]
=(e^{\al} D_{ij}^{(k)}l)(f)-
l(f(h+\al(h))(\eii^{(k)}-\ejj^{(k)}))
\\[10pt]
=(e^{\al}D_{ij}- e^{\al}(
\pp_i^{(k)}- \pp_j^{(k)}))(l)(f).
\end{array}
$$ 
This proves i).

iii) We have
$$
\renewcommand{\arraystretch}{1.4}
\begin{array}{l}
\fii*(f\cdot \fii^{(k)})=
\fii f\cdot \fii^{(k)}- f\cdot \fii^{(k)}\fii
\\[10pt]
=f(\fii \fii^{(k)}- \fii^{(k)}\fii)
=f(2 \fii \fii^{(k)}-[\fii, \fii^{(k)}])
=2f \fii\cdot \fii^{(k)}-2\eii^{(k)}f
\end{array}
$$ 
or $f\cdot \fii\cdot \fii^{(k)}=\ds\frac12\fii*( f\cdot \fii^{(k)})+
f\cdot \eii^{(k)}$.
This proves iii). Heading iv) is similar. \end{proof}

\ssbegin{3.8.2}{Corollary} Heading {\em iii)} of Theorem $1.5$ is true.
\end{Corollary}

\ssec{3.9.  Proof of heading iv) of Theorem 1.5} Let
$V^{\lambda}=\Span (v_p)$, where the $v_p$ form a basis of
$V^{\lambda}$, let the $v^*_p$ form the left dual basis; set
$\omega_{\lambda}^*=\sum v_p^*\ot v_p$, and $\omega_{\lambda}=\sum v_p\ot
v_p^*$.

If $\fg$ is embedded into $\fg\oplus \fg$ as the first summand, and
$u\in U(\fg)$, then
$$
\renewcommand{\arraystretch}{1.4}
\begin{array}{l}
\vh_{\lambda}(u\ot 1)=\ttt(\omega_{\lambda}^*, \omega_{\lambda})(u\ot
1) =\omega_{\lambda}^*(u\ot 1(\omega_{\lambda})) \\[10pt]
=\omega_{\lambda}^*(\ds\sum uv_p\ot v_p^*) =\ds\sum v_p^*(uv_p)=tr(u).
\end{array}
$$ 
In other words, the matrix coefficient
$\ttt^{\pi}(\omega_{\lambda}^*, \omega_{\lambda})$ coincides with the
functional $\tr_{V^{\lambda}}(u)$ after restriction onto $U(\fg)$. But
due to \cite{Se1} we have, up to a scalar multiple, 
$$
\tr_{V^{\lambda}}|_{\fh_{\bar 0}}=Q_{\lambda}(e^{\vr_1}, \ldots ,
e^{\vr_n}).\qed
$$ 

\ssec{3.10.  Proof of headings ii) and iii) of Lemma 1.2.3} Heading
ii) immediately follows from our proof of Corollary 3.4.  Let us prove
iii).  It is not difficult to see that the algebra $\Omega$ generated
by all the $\Omega_k$, $k=1, 3, 5, \ldots$, is the image of $Z(\fq(n))$
under the homomorphism
$$
r: Z(\fq(n))\tto\Omega,\quad r(z)Q_{\lambda}=\vh(z)(
\lambda)Q_{\lambda}, 
$$ 
where $\vh$ is the Harish--Chandra homomorphism, see \cite{Se3}. 
Therefore, statement of heading iii) can be reformulated as follows:

{\sl Let $R_n$ be the algebra of polynomials $r(t_1, \ldots , t_n)$
symmetric with respect to $(t_1, \ldots , t_n)$ and which do not
depend on $s$ after substitution $t_i=s$, $t_j=-s$. Define an
$R_{n}$-action on the algebra generated by $Q_{\lambda}$ by setting
$$
r*Q_{\lambda}=r(\lambda)Q_{\lambda}.
$$ 
If $P$ is an eigenvector for $R_n$ and belongs to the subalgebra
generated by $Q_{\lambda}$, then $P$ coincides, up to a scalar
multiple, with one of the $Q_{\lambda}$.}

Indeed, let $P=\sum_{\lambda\in \Lambda}C_{\lambda }Q_{\lambda }$,
where $\Lambda$ is the set of partitions of length $\leq n$.  Since
$r*P=C(r)P$ it follows that $r(\lambda)=const$ for any $\lambda \in
\Lambda$.  Let $\lambda, \mu\in\Lambda $ and $\lambda \ne \mu$; then
$r({\lambda })=r(\mu)$ for any $r\in R_{n}$.  Hence,
$$
\prod\frac{t-\lambda_i }{t+\lambda_i }=\prod\frac{t-\mu_j}{t+\mu_j}. 
$$ 
We may assume that $\lambda_i>0 $ and $\mu_i>0$ for all $i$ and $j$. 
Then the identity 
$$
\prod(t-\lambda_i )\prod(t+\mu_j)= \prod(t+\lambda_i)\prod(t+\mu_j)
$$ 
implies $\prod(t-\lambda_i)=\prod(t-\mu_j)$, hence, $\lambda=\mu$. 
Contradiction. \qed

\section*{\S 4. Bispherical functions and the radial parts of Laplace 
operators for the triple $(\fpe(n), \fgl(n|n), \fpe(n))$}

\ssec{4.1} Let $V$ be an $(n|n)$-dimensional superspace, $I=\Io\cup
\IIo=\{1, \ldots , n\}\cup\{\bar 1, \ldots , \bar n\}$, the union of
the even and odd indices, see sec. 1.4.1. Let $\{e_{i}\}_{i\in I}$
be the basis of $V$ consisting of vectors whose parity is equal to
that of their indices and $\{e_{ij}\}_{i, j\in I}$ be the basis fo
$\Mat (V)$ consisting of matrix units, cf. sec. 1.4.2. 

The supertransposition antiautomorphism in $\fgl(V)$ is in these terms
of the form
$$
\eij^t=(-1)^{p(i)(p(j)+\bar 1)}\eji. \eqno(4.1.1)
$$ 
Define two 
antiautomorphisms $\psi_1$ and $\psi_2$ by setting
$$
\psi_1(x)=(-1)^x\Pi(x^t)\Pi\; \text{ and } \;\psi_2=P\psi_1P.\eqno(4.1.2)
$$ 
Direct calculations show that
$$
\psi_1(\eij)=(-1)^{p(j)(p(i)+\bar 1)}e_{\bar \jmath\bar i}, \quad 
\psi_2(\eij)=(-1)^{p(i)(p(j)+\bar 1)}e_{\bar \jmath\bar i}.
\eqno(4.1.3)
$$ 
Define two Lie subsuperalgebras of $\fgl(V)$:
$$
\fpe_1(V)=\{x\in \fgl(V)\mid \psi_1(x)=-x\}\eqno(4.1.4)
$$ 
and 
$$
\fpe_2(V)=\{x\in \fgl(V)\mid \psi_2(x)=-x\}.\eqno(4.1.5)
$$ 
Observe that $\psi_1(x)=\psi_2(x)$ if $p(x)=0$ and $\psi_1(x)=-\psi_2(x)$ 
if $p(x)=1$. Therefore, 
$$
\fpe_1(V)=\fgl(V)_{\bar 0}^-\oplus \fgl(V)_{\bar 1}^-, \quad 
\fpe_2(V)=\fgl(V)_{\bar 0}^-\oplus \fgl(V)_{\bar 1}^+, 
$$ 
where
$$
\renewcommand{\arraystretch}{1.4}
\begin{array}{l}
 \fgl(V)_{\bar 0}^-=\{x\in \fgl(V)_{\bar 0}\mid \ \psi_2(x)=-x\}, \\
\fgl(V)_{\bar 1}^-=\{x\in \fgl(V)_{\bar 1}\mid \ \psi_1(x)=-x\}, \\
\fgl(V)_{\bar 1}^+=\{x\in \fgl(V)_{\bar 1}\mid \ \psi_1(x)=x\}. 
\end{array}
$$
For every $x\in\fgl(V)$, set
$$
x^+=\frac12(x+\psi_1(x)), \quad x^-=\frac12(x-\psi_1(x)).
$$ 
Also set
$$
\fh^+=\Span(\eii^+\mid i\in I_{\bar 0}).
$$

\ssbegin{4.2}{Lemma} For $f\in S(\fh^+)$ and $\alpha
=\eps_{i}-\eps_{j}$ set
$$
\renewcommand{\arraystretch}{1.4}
\begin{array}{l}
R_{ij}^-f=\frac12[f(h-\al(h))-f(h+\al(h))], \\
R_{ij}^+f=\frac12[f(h-\al(h))+f(h+\al(h))].\end{array}
$$ 
Then the following identities hold:

{\em i)} $\eij^-f=R_{ij}^+ f\eij^-+R_{ij}^- f\cdot\eij^+$.

{\em ii)} $R_{ij}^- f\eij\eji^{(k)}-(R_{ij}^--R_{ij}^+) f\cdot
[\eij^-, \eji^{(k)}]\in \fpe_1U(\fg)+U(\fg)\fpe_2$ if $p(i)+p(j)=0$.

{\em iii)} $R_{ij}^+ f\eij\eji^{(k)}-(R_{ij}^+-R_{ij}^-) f\cdot
[\eij^+, \eji^{(k)}] \in \fpe_1 U(\fg)+U(\fg)\fpe_2$ if $p(i)+p(j)=1$.

{\em iv)} $[\eij^-, \eji^{(k)}]=\ds\frac12[\eii^{(k)}-\ejj^{(k)}]$ if
$p(i)+p(j)=0$.

{\em v)} $[\eij^+, \eji^{(k)}]=\ds\frac12[\eii^{(k)}+\ejj^{(k)}]$ if
$\bar i\ne j$ and $p(i)+p(j)=1$.

{\em vi)} $[e_{i\bar i}^{+}, \eii^{(k)}]=\eii^{(k)}+e_{\bar
i\bar i}^{(k)}$ if $p(i)=1$. 

{\em vii)} $[e_{i\bar i}^{+}, e_{\bar ii}^{(k)}]=0$ if $p(i)=0$.
\end{Lemma}

Proof is reduced to a direct verification. \qed

\ssec{4.3. Proof of heading i) of Theorem 1.7} It suffices to prove 
that $U(\fg)=S(\fh^+)+\fpe_1U(\fg)+U(\fg)\fpe_2$.

Indeed, any element of $U(\fg)$ can be represented in the form 
$$
u=u_0+u_1u'_0u_2, \text{ where } u_1\in U(\fg_{\bar 1}^-), \ 
u_2\in U(\fg_{\bar 1}^+), \ u_0, u_0'\in U(\fg_{\bar 0}). 
$$ 
Therefore, we may assume that 
$$
u=f e^+_{\al_1} \ldots e^+_{\al_k}, 
$$ 
where $f\in S(\fh^+)$, and $e^+_{\al_1}, \ldots , e^+_{\al_k}\in \fg_{\bar
0}$ are the weight vectors.

Now, induction on $k$.  If $k=1$, then $R_{ij}^- f \eij^+=\eij^-
f-R_{ij}^+ f \eij^-$ by Lemma 4.2.i).  Hence, $f \eij^+\in \fg^-_{\bar
0}U(\fg_{\bar 0})\fg^-_{\bar 0}$.  Let $k>1$.  Then
$$
\renewcommand{\arraystretch}{1.4}
\begin{array}{rcl}
R^-_{{\al_1} }fe^+_{\al_1} \ldots e^+_{\al_k} &=&
e^-_{\al_1}fe^+_{\al_2}\ldots e^+_{\al_n}-R^+_{{\al_1} }
fe^-_{\al_1}e^+_{\al_2}\ldots \\
[10pt] &\equiv_{\pmod {\fg^-_{\bar
0}U(\fg)\fg^-_{\bar 0}}}& -R^+_{{\al_1} }f e^-_{\al_1}e^+_{\al_2}
\ldots e^+_{\al_n} \\ 
[10pt] &\equiv& R^+_{\al_1} f\cdot [e^-_{\al_1}, 
e^+_{\al_2} \ldots e^+_{\al_k}] \\ 
[10pt] &=& -R^+_{\al_1}
f\cdot[e^-_{\al_1}, e^+_{\al_2}] e^+_{\al_3}\ldots e^+_{\al_k} \\
[10pt]
&& -R^+_{\al_1} f e^+_{\al_2}[e^-_{\al_1}, e^+_{\al_3}]\ldots
e^+_{\al_n}+\ldots \in \fg^-_{\bar 0} U(\fg_{\bar 0}) \fg^-_{\bar 0} .
\end{array}
$$

\ssec{4.4. Proof of Lemma 1.2.1} Set
$$
\dekip=\tilde\pp_i^{(k)}+ \tilde\pp_{\bar i}^{(k)}, \quad\quad 
\dekim=\tilde\pp_i^{(k)}- \tilde\pp_{\bar i}^{(k)}.
$$ 
Now, it is not difficult to verify the following identities:
$$
\dekip-\dekiip=\pp_i\dekiim+\sum_{j\ne i}\frac{2}{e^{\vr_{ij}}- 
e^{\vr_{ji}}}(\dekiim-\dekjim), \eqno(4.4.1)
$$
$$
\dekim=(\pp_i-1)\dekiip +\sum\left(\frac{2}{e^{\vr_{ij}}- 
e^{\vr_{ji}}}\dekiip -\frac{2 e^{\vr_{ij}}}
{e^{\vr_{ij}}- e^{\vr_{ji}}}\dekjip\right).\eqno(4.4.2)
$$
Further, if $\Delta_i^{(2k)+}=f_{2k}(\pp_i^{(1)}, \ldots ,
\pp_i^{(2k-1)})$ and $\Delta_i^{(2k+1)^+}=f_{2k+1} (\pp_i^{(1)}, \ldots
, \pp_i^{(2k+1)})$, where $f_{2k}$ and $f_{2k+1}$ are any linear
functions, then
$$
\renewcommand{\arraystretch}{1.4}
\begin{array}{rcl}
f_{2k+1}&=& f_{2k}+ f_{2k-1}(\pp_i^{(3)}, \ldots , \pp_i^{(2k+1)}), \\
f_{2k}&=& f_{2k-1}+ f_{2k-2}(\pp_i^{(3)}, \ldots , \pp_i^{(2k-1)}).\\
\end{array}
\eqno(4.4.3)
$$
It is easy to verify that 
$$
\tilde\Omega_1=2\Omega_1=\tilde\Omega_2, \quad
\tilde\Omega_3=2\Omega_3+2\Omega_1, \quad
\tilde\Omega_4=4\Omega_3+2\Omega_1
$$ 
and, by induction, (4.4.3) implies that $\tilde\Omega_k$ is a linear
combination of the $\Omega_{2l+1}$ for $2l+1\leq k$.  One can show
that, the other way round, $\Omega_{2l+1}$ can be expressed via
$\tilde\Omega_k$.  \qed

\ssec{4.5.  Proof of headings ii), iii) of Theorem 1.5} Statement of
heading ii) is obvious.  Heading iii) follows from Lemma 4.2 and the
fact that $z_k=\sum_{i\in I}\eii^{(k)}$.

\ssec{4.6.  Proof of heading iv) of Theorem 1.5} It is easy to verify
that $\ttt^*=\sum e_i^*\ot e_{\bar i}^*$ is a $\fpe_1$-invariant
whereas $\ttt=\mathop{\sum}\limits_{i\in I}e_i\ot e_{\bar i}$ is a
$\fpe_2$-invariant.  According to \cite{Se2} the linear hull of all
the $\fpe_1$-invariants in $V^{\ot 2k}$ is isomorphic to
$$
\Ind_{H_k}^{\fS_{2k}}(\vr)=\mathop{\oplus}\limits_{\lambda=(\al_1+1,
\ldots , \al_p+1|\al_1, \ldots , \al_p)} S^{\lambda}, \eqno(4.6.1)
$$
where $H_k={\fS}_k\circ \Zee_2^k$ is the semidirect product, and
$h((\ttt^*)^{\ot k})=\vr(h)(\ttt^*)^{\ot k}$ for any $h\in H_k$.

Similarly, the module of all the $\fpe_2$-invariants in $V^{\ot 2k}$ is 
of the form (4.6.1). This implies that
$$
\vh_{\lambda}(u)=\ttt(v_{\lambda}^*, \;
v_{\lambda})(u)=\ttt(e_{\lambda}\ttt^{*2k},\; e_{\lambda}\ttt^{2k})(u), 
$$ 
where $e_{\lambda}$ is the minimal idempotent in the Hecke algebra
$H(\fS_{2k}, H_{2k}, \vr)$ corresponding to partition $\lambda$, see
\cite{St}.

If $\sigma\in \fS_{2k}$, then the map
$$
\sigma\mapsto \vh_{\sigma}, \ \text{ where }\ 
\vh_{\sigma}(u)=(-1)^{\frac12k(k-1)}\ttt^{*\ot k}(\sigma u\ttt^{\ot
k}),
$$ 
satisfies
$$
\vh_{h_1\sigma h_2}=\vr(h_1h_2)\vh_{\sigma}\ \text{ for any
$h_1, h_2\in H_k$. }
$$
Therefore, we obtain a map 
$H(\fS_{2k}, H_k, \vr)\tto U(\fg)^{*inv}$. By restricting this map
onto $S(\fh_{\bar 0})$ we obtain a map
$$
\ch: H(\fS_{2k}, H_k, \vr)\tto S(\fh_{\bar 0})^*.\eqno(4.6.2)
$$ 
Let $K=K_1\cup K_2$, where $K_1=\{1, \ldots , 2k_1\}$, $\ K_2=\{
2k_1+1, \ldots , 2k_1+2k_2\}$ and let $\sigma_1$ permute the elements of
$K_1$ whereas $\sigma_2$ permutes the elements of$K_2$; let $\sigma
=\sigma_{1}\sigma_{2}$. Now, one can verify that
$$
\ch(\sigma_{1}\sigma_{2})=\ch(\sigma_{1})\ch(\sigma_{2}).
$$ 
Let now $k$ be odd and $\sigma =(1, \ldots , 2k)$ a cycle. Let us
calculate $\ch(\sigma )$. We have
$$
\ttt^{\ot k}=\sum e_{\psi}, \quad \ttt^{*\ot k}=\sum e_{\psi}^*, \eqno(4.6.3)
$$ 
where the sums run over all the maps $\psi:\{1, \ldots , 2k\}\mapsto
\{1, \ldots , n, \bar 1, \ldots , \bar n\}$ such that
$\psi(2i)=\overline {\psi(2i-1)}$ for $i=1, \dots , k$ and where
$e_{\psi}=e_{\psi_{(1)}}\ot\ldots\ot e_{\psi_{(2k)}}$ and $e^*_{\psi}$
is similarly defined.  Therefore,
$$
\ch(\sigma)(u)=(-1)^{\frac12k(k-1)}\ttt^{*\ot k}(u\sum \sigma e_{\psi}), 
$$ 
where in the sum one has to take into account only the summands
$e_{\psi}$ for which $\sigma \psi$ possesses the same property as
$\psi$ does, i.e., $\psi(2i)=\overline{\psi(2i-1)}$.  Therefore, we
may assume that
$$
\psi(1)=\overline {\psi(2)}= \psi(3)=\overline {\psi(4)}=\ldots=
\psi(2k-1)=\overline{\psi(2k)}.
$$ 
i.e., $e_{\psi}=(e_i\ot e_{\bar i})^{\ot k}$ for $i\in
\{1, \ldots , n, \bar 1, \ldots , \bar n\}$. But for such $\psi$ and $u\in
S(\fh_{\bar 0})$ we have
$$
\renewcommand{\arraystretch}{1.4}
\begin{array}{rcl}
ch(\sigma)(u)
&=&
(-1)^{\frac12k(k-1)}\ttt^{*\ot k}(u\sum_{i\in I}\sigma(e_i\ot e_{\bar i})^{\ot k})
\\[10pt]
&=&
(-1)^{\frac12k(k-1)}\ttt^{*\ot k}(u\sum_{i\in I}(e_i\ot e_{\bar i})^{\ot k})
\\[10pt]
&=&
(-1)^{\frac12k(k-1)}\ttt^{*\ot k}(\sum_{i\in I_{\bar 0}}2e^{\vr_i}(u)
(e_i\ot e_{\bar i})^{\ot k})
\\[10pt]
&=&
2\sum_{i\in I_{\bar 0}}2e^{k\vr_i}(u).
\end{array}
$$ 
Hence, $\ch(\sigma)= 2\sum_{i=1}^n2e^{k\vr_i}$. Therefore, if 
$\sigma=\sigma_{\nu_1} \ldots \sigma_{\nu_p}$ is the product of
independent cycles of odd lengths, then
$$
\ch(\sigma)=2^{l(\nu)}P_{\nu}, \text{ where } 
P_{\nu}= P_{\nu_1}\ldots P_{\nu_p} \text{ and } P_l=\sum_{i=1}^ne^{l\vr_i}.
$$ 
We see that the map (4.6.2) coincides with the characteristic map
Stembridge constructed in \cite{St}, p. 85. Therefore, by Theorem
5.2 from \cite{St} we have
$$
\ch(e_{\lambda})=\ttt_{\lambda}(e^{\vr_1}, \ldots , e^{\vr_n})\cdot
2^{n-e(\lambda)} \cdot g_{\lambda}, 
$$ 
where $g_{\lambda}$ is the number of shifted standard tableaux of
shape $\lambda$. \qed

\end{document}